\documentclass[12pt,reqno]{amsart}
\newcommand{\ie}{{\em i.e.} }
\newcommand{\eg}{{\em e.g.} }
\newcommand{\EV}{eigenvalue }
\newcommand{\EVP}{eigenvalue problem }
\newcommand{\EVs}{eigenvalues }

\newcommand{\EFs}{eigenfunctions }

\newcommand{\sa}{self-adjoint }
\newcommand{\nsa}{non-self-adjoint }
\newcommand{\BC}{boundary condition }
\newcommand{\DBC}{Dirichlet boundary conditions }
\newcommand{\NBC}{Neumann boundary conditions }
\newcommand{\BCs}{boundary conditions }
\newcommand{\BVP}{boundary value problem }

\newcommand{\Real}{\mathbb R}

\newcommand{\Lom}{\mathcal{L}}

\def\R{{\Bbb R}}
\def\C{{\Bbb C}}

\def\rmd{{\rm d}}
\def\rme{{\rm e}}
\def\dpa{\partial}
\def\re{{\rm Re}}
\def\im{{\rm Im}}

\def\Arg{{\rm Arg}}

\def\a{\alpha}
\def\b{\beta}
\def\g{\gamma}
\def\d{\delta}
\def\e{\varepsilon}

\def\y{\eta}

\def\l{\lambda}

\def\n{\nu}
\def\x{\xi}

\def\f{\varphi}

\def\o{\omega}

\def\G{\Gamma}
\def\D{\Delta}

\def\O{\Omega}

\theoremstyle{definition}

\theoremstyle{remark}

\numberwithin{equation}{section}

\begin{document}

\title[Perturbed cylinders]{Separation of variables in
perturbed cylinders}
\author{A. Aslanyan and E.B. Davies}
\address{
Department of Mathematics\\
King's College London\\
Strand, London WC2R 2LS, UK}
\email{aslanyan@mth.kcl.ac.uk \\
e.brian.davies@kcl.ac.uk}
\thanks{We thank Engineering
and Physical Sciences Research Council for support under grant No.
GR/L75443.}
\date{11 December 2000}


\begin{abstract}
We study the Laplace operator subject to
\DBC in a two-dimensional
domain that is one-to-one mapped onto a cylinder
(rectangle or infinite strip). As a result of this
transformation the original \EVP is reduced to an
equivalent problem for an operator with variable
coefficients. Taking advantage of the simple geometry
 we separate variables by means of the
Fourier decomposition method. The ODE system obtained
in this way is then solved numerically yielding the
\EVs of the operator. The same approach
allows us to find complex resonances arising in some
non-compact domains. We discuss numerical examples related to
quantum waveguide problems.

\end{abstract}

\maketitle
\section{Introduction}

The object of this study is the Dirichlet Laplacian in a perturbed cylinder, \ie
a domain that is mapped onto a rectangle or an
infinitely long strip depending on the domain being
compact or non-compact. A typical example of a perturbed cylinder
is a waveguide where the propagation of waves
is governed by the Helmholtz equation. The two major  types of
waves observed in waveguides are
referred to as trapped modes and resonance solutions.
Mathematically trapped modes are described as the \EFs
of \sa operators. These are $\Lom^2$ functions associated
with bound states, as opposed to resonance solutions
that correspond to the so-called scattering poles, or
complex resonances. Despite their different nature,
\EVs and resonances are sometimes closely connected with
each other. For example, \EVs may generate resonances
under small perturbations of the domain. It is this
situation that interests us and motivates our study
of resonances.

Waveguide phenomena are usually
associated with either Dirichlet or Neumann boundary
conditions. The former correspond to scattering problems in quantum theory,
the latter appear in acoustics. We refer to the series of papers \cite{Ex1,Ex2,Ex3} where
bound states and scattering in quantum waveguides are studied.
In \cite{SRMPP} the authors address the problem of finding quantum resonances
numerically for a particular waveguide. For results on acoustic waveguides, both
theoretical and numerical see, for instance,
\cite{EvLeVa,MaNe,DaPa,APV,AbKr,OtLa}. The papers cited here are
concerned with either \EVs or resonances occurring in waveguides
under specific conditions, or both of these. Our paper
is very close in spirit to \cite{APV} and \cite{AbKr}. In the former
the main issue is the resonance--\EV connection, and the
technique of the latter is also based
on the separation of variables.

Our intention is to study the above mentioned  problems
numerically. Dealing with both of them involves solving
a \BVP for the Dirichlet Laplacian in two dimensions,
that is either a \sa \EVP or a \nsa resonance problem.
Along with general numerical methods applicable in two dimensions,
there exist techniques especially designed for cylinder-like
domains also called ducts in acoustics. We have already mentioned
\cite{AbKr,OtLa} where such methods are developed.
Both of these papers stress the importance of
advanced methods specially designed for acoustic waveguides.
It is hardly surpsising that carefully performed numerical analysis
is equally important for quantum problems. In \cite{OtLa} the authors
apply a second-order finite difference method and implement an
iterative procedure for the resulting algebraic system.
It is mentioned there that standard methods not using any
preconditioning are likely to fail especially when a large wave
number is involved.

The numerical approach
proposed in \cite{AbKr} is similar in spirit to that of our paper.
In both cases the Helmholtz equation is reduced to
the so-called coupled mode system of equations via the separation of variables. The main
difference is that in \cite{AbKr} the coefficients of the
ODE system have to be computed numerically whereas our
choice of the Fourier expansion functions allows us to find
them in closed form. The Dirichlet problem studied here is separable as
opposed to the  much less straightforward Robin case.
The examples in the cited paper
are related to higher frequencies while we concentrate
on the lowest oscillation mode only. On the other hand, the
transfer method we use for the final ODE problem is able to
handle a waveguide with a narrow throat --- a situation
not covered in \cite{AbKr}.

The aim of this paper is to elaborate a method suitable for
perturbed cylinders that takes account of their
geometry. The method based on the Fourier decomposition
in one direction allows us to separate variables in the
Helmholtz equation explicitly leading to a system of  ODEs. This
is done for a fairly general geometry in the next section.
In section 3 we discuss  different \BCs involved. First we deal with standard \sa conditions,
then concentrate on non-compact domains and define resonances by a
specific \BC at infinity. What is often called the radiation condition in the
literature is  rewritten in terms of the resulting ODE
system. We end up with a \nsa \EVP on a finite interval
whose solution  approximates that of the original
resonance problem. Finally, we use the transfer method
of \cite{Abr1} to find the \EVs of the two
problems. Our numerical results illustrate the closeness of
\EVs and resonances and are presented in section 4 where we
conclude by discussing the rate of convergence.

To be able to compare the method of the paper with others we look at the finite volume
method described, \eg in \cite{Var}. As discussed in section \ref{ne},
the proposed approach tested on our \EV examples proves to be significantly
more efficient than the standard two-dimensional procedure.

\section{Separation of variables for the Laplacian}

\subsection{Change of variables}\label{sb1}

Consider the operator $H:=-\D$ acting on $\Lom^2(\O)$
subject to Dirichlet boundary conditions. The domain
$\O$ is defined as
\begin{equation}
\O = \{(\x,\y): \ a<\x<b, \ 0<\y<\f(\x)\} \label{dom}
\end{equation}
in the Cartesian coordinates $(\x,\y)$. The possibility of
$a$ and $b$ being infinite is not excluded here, so that
$\O$ is not necessarily compact. The function $\f(\x)$ is
assumed to be smooth and satisfy $\f(\x)>0, \
\x\in[a,b]$. To find the spectrum of $H$ we solve the
Helmholtz equation
\begin{equation}
-\D f(\x,\y)  = \l f(\x,\y), \qquad (\x,\y) \in \O \label{hel}
\end{equation}
with the \BCs
\begin{equation}
f(\x,\y)\,=\,0\,, \qquad (\x,\y) \in \partial\O \label{dbc}
\end{equation}

The change of variables
\begin{equation}
x=\x, \qquad y=\y/\f(\x) \label{cv}
\end{equation}
maps the perturbed cylinder $\O$ onto $\O_0\,=\,\{(x,y): \ a < x < b, \ 0< y <1\}$
which is either a rectangle or a strip (infinite or semi-infinite).
We mention that a similar method has been used by Borisov et al. \cite{BEGK}
to study bound states associated with a local perturbation of a strip or
layer. In our case the deformation reduces the width of the strip locally,
and there are no bound states. The transformation (\ref{cv}) can be expressed in the differential form as
$$
\nabla_{\x\y} = W \nabla_{xy},
\qquad \nabla_{xy}=
\left( \begin{array}{c}
\frac{\dpa}{\dpa x} \\
\frac{\dpa}{\dpa y}
\end{array} \right),
\qquad W=
\left( \begin{array}{cc}
1 & -\frac{\f'y}{\f} \\
0 & \frac{1}{\f}
\end{array} \right).
$$
The quadratic form corresponding to $H$ is given by
$$
J(f) \,=\, \int_{\O} \left[ (\nabla_{\x\y}f,\nabla_{\x\y}f) - \l(f,f)
\right] \rmd\x\rmd\y
$$
or, equivalently, by
$$
J(f) \,=\,\int_{\O_0} \left[ (W\nabla_{xy}f,W\nabla_{xy}f) - \l(f,f)
\right] \f(x)\rmd x \rmd y.
$$
This can be rewritten as
\begin{equation}
J(f)=\,\int_{\O_0} \left[ (A\nabla_{xy}f,\nabla_{xy}f) -
\l(f,f)
\right] \f(x)\rmd x \rmd y,
\label{qf}
\end{equation}
where
\begin{equation}
A\,=\,W^*W\,=\,
\left( \begin{array}{cc}
1 & -\frac{\f'y}{\f} \\
-\frac{\f'y}{\f} & \frac{(1+(\f'y)^2)}{\f^2}
\end{array} \right). \label{a}
\end{equation}
Hence the Helmholtz equation in the new variables takes the form
\begin{equation}
\begin{array}{r}
{\displaystyle \left(\f\, f_x\right)_x\,
-\, \left(\f'y\,f_x\right)_y
\,-\, \left(\f'y \,f_y\right)_x\,}
{\displaystyle +\,
\left((1+(\f'y)^2)\,f_y/\f\right)_y\, +\,
\l \f f \,=\, 0.}
\end{array}
\label{ee}
\end{equation}

Since we restrict ourselves to the Dirichlet case no change
in the \BCs is required here, the condition (\ref{dbc}) is
retained on $\partial\O_0$. However, in a
generic situation one can still use (\ref{cv}), (\ref{ee}) provided
obvious changes are made to the original \BCs where
necessary. For instance, instead of \NBC at $\y=\f(\x)$ one would have
$$
\left. \left( (1+\f^{\prime 2})f_y\,-\,\f\f'\,
f_x \right) \right|_{y=1}\,=\,0,
$$
while \NBC at $x=a,\,b$ would become
$$
\left. (\f f_x - \f' y f_y)\right|_{x=a,\,b}
$$

Note that a similar transformation can be also done for a more
general domain
$$
\O = \{(\x,\y): \ a<\x<b, \ \f_1(\x)<\y<\f_2(\x)\}.
$$
The change of variables
$$
x=\x, \qquad y=\frac{\y-\f_1}{\f_2-\f_1}
$$
leads to a quadratic form of type (\ref{qf}) whose coefficients are not quoted here for brevity.

\subsection{Discretisation in the $y$-direction}\label{sb2}

The quadratic form (\ref{qf}) is related to the
transformed operator on the weighted space
${\Lom}^2(\O_0,\f\rmd x \ \rmd y)$. Here and below
we use the  notation $\nabla=\nabla_{xy}$.

To discretise the form $J(f)$ in the $y$-direction let us
separate the variables expanding $f$ as
\begin{equation}
f(x,y) \,=\, \sum_{k=1}^\infty\,  g_k(y) h_k(x). \label{ex}
\end{equation}
Recall that we have \DBC everywhere so that our
natural choice is to work with an orthonormal system
of functions vanishing at the horizontal parts of the
boundary. We therefore opt for
$$
g_k\,=\,\sqrt{2}\sin(\pi k y)\,.
$$
Denote $f_0=f, \ f_1=f_x, \ f_2 = f_y$;
$$
f_i = \sum_{k=1}^\infty\, h_k^i(x) g_k^i(y)\,,
$$
where
$$
g_k^0\,=\,g_k^1\,=\,g_k, \qquad
g_k^2\,=\,\sqrt{2}\cos(\pi k y)\,.
$$
In this notation the $x$-dependence is determined by
the functions
$$
h_k^0 = h_k, \qquad h_k^1 = h'_k, \qquad h_k^2 = \pi k
h_k
$$
(throughout the paper $'$ denotes differentiation with respect to
$x$).

We notice that the variables are separated in the coefficients
of $J(f)$:
\begin{equation}
\f(x)A_{ij}(x,y) = B_{ij}(x) + C_{ij}(y)D_{ij}(x), \qquad i,j=1,2.
\label{dec}
\end{equation}
The entries $A_{ij}$ of the matrix $A$ are defined
by (\ref{a}); the matrices $B, \ C$ and $D$ satisfying
the above decomposition are given below:
$$
B\,=\,
\left( \begin{array}{cc}
\f & 0 \\
0 & \frac{1}{\f}
\end{array} \right), \qquad
C\,=\,
\left( \begin{array}{cc}
0 & y \\
y & y^2
\end{array} \right), \qquad
D\,=\,
\left( \begin{array}{cc}
0 & -\f'\\
-\f' & \frac{\f^{\prime 2}}{\f}
\end{array} \right).
$$

The formula (\ref{dec}) allows us to rewrite (\ref{qf}) in the form
$$
J(f)\,=\,\int_a^b\int_0^1 \left[\sum_{i,j=1}^2\, (B_{ij}
+ C_{ij}D_{ij}) \bar{f}_i f_j - \l \f |f|^2
\right] \rmd x \rmd y
$$
$$
=\,\int_a^b \left[\sum_{i,j=1}^2\,B_{ij} \int_0^1 \bar{f}_i f_j
\rmd y
+ D_{ij} \int_0^1 C_{ij} \bar{f}_i f_j \rmd y - \l \f \int_0^1 |f|^2 \rmd y \right]\rmd x
$$
$$
=\,\int_a^b \left[\sum_{i,j=1}^2\,B_{ij} E_{ij}
+ D_{ij} F_{ij}  - \l \f G \right]\rmd x.
$$
As we substitute the expansion (\ref{ex}) into the
above integral, the coefficients $E_{ij}, \ F_{ij}, \ G$ are
readily computed below.
$$
E_{ij}\,=\,\int_0^1 \bar{f}_i f_j \rmd y \,=\,\int_0^1
\sum_k \bar{h}_k^i g_k^i \sum_r h_r^j g_r^j \rmd y \,
$$
$$
=\,\sum_{k,r} \bar{h}_k^i h_r^j \int_0^1 g_k^ig_r^j \rmd y\,=\,
\sum_{k,r}  \a_{kr}^{ij} \bar{h}_k^i h_r^j;
$$
$$
F_{ij}\,=\,\int_0^1 C_{ij}\bar{f}_i f_j \rmd y \,=\,
\int_0^1 C_{ij} \sum_k \bar{h}_k^i g_k^i \sum_r h_r^j
g_r^j \rmd y
$$
$$
=\,\sum_{k,r} \bar{h}_k^i h_r^j \int_0^1 C_{ij}g_k^ig_r^j \rmd y\,=
\,\sum_{k,r}  \b_{kr}^{ij} \bar{h}_k^i h_r^j\,;
$$
$$
G\,=\,\int_0^1 \sum_k \bar{h}_k^0 g_k^0 \sum_r h_r^0
g_r^0 \rmd y \,=\, \sum_k |h_k|^2\,.
$$
To find $\a_{kr}^{ij}$ we use the orthogonality
relations for $g_k^i$; in fact, we only need the
diagonal elements $\a_{kr}^{jj}=\d_{kr}$. The coefficients
$\b_{kr}^{ij}$ are also calculated in the closed form:
$$
\b_{kr}^{11}\,=\,0, \qquad \b_{kr}^{22}\,=\,
\begin{cases}
\frac{1}{3} + \frac{1}{2\pi^2k^2}, & k=r \\
(-1)^{k+r}\frac{4(k^2+r^2)}{\pi^2(k^2-r^2)^2}, & k\not=r
\end{cases}\,,
$$
$$
\b_{kr}^{12}\,=\,\b_{rk}^{21}\,=\,
\begin{cases}
 - \frac{1}{2\pi k}, & k=r  \\
(-1)^{k+r}\frac{2k}{\pi(r^2-k^2)}, & k\not=r
\end{cases}\,.
$$
The quadratic form is now reduced to that of a one-dimensional
differential problem.

\subsection{Canonical ODE system}
Having done the above calculations we finally arrive at
$$
J(f)\,=\,\int_a^b \left[ \sum_k \,|h'_k|^2 \,+\, \left(
\left[\frac{\pi k}{\f}\right]^2 - \l\right)|h_k|^2
\right.
$$
$$
\left. +\, \sum_{k,r} \,\frac{\f^{\prime 2}}{\f^2}\pi^2kr \b_{kr}^{22}
\bar{h}_kh_r \,-\, \frac{\f'}{\f}\pi \left( r \b_{kr}^{12}
\bar{h}'_kh_r + k \b_{rk}^{12} \bar{h}_kh'_r \right)
\right] \f\, \rmd x. \label{qfx}
$$
The Euler equations are easily derived in the
standard way. A simple calculation shows that the discretised
form (\ref{qfx}) is equivalent to the ODE system written in
its canonical self-adjoint form as
\begin{equation}
-(Ph')' \,+\, Qh' \,-\, (Q^*h)' \,+\, Rh\,=\,0. \label{sl}
\end{equation}
Here the vector of unknowns
$$
h\,=\,(h_1, h_2, \ldots)^T;
$$
the matrix coefficients are given by
$$
\qquad P_{kr}\,=\,\f \d_{kr},
\qquad Q_{kr}\,=\, -\pi k \b_{rk}^{12}\f',
$$
$$
R_{kr}\,=\,\pi^2 kr \b_{kr}^{22}\frac{\f^{\prime 2}}{\f} \,+\,\left( \frac{(\pi k)^2}{\f} -
\l\f\right)\d_{kr}, \qquad k,r=1,2,\ldots \,.
$$
For practical purposes we truncate the system to a
finite number of equations taking a sufficiently large
$N$ and keeping the same notation $h, \ P, \ Q, \ R$
for the truncated  matrices where $k,r=1,\ldots,N$.
This is justified by the fact that the Fourier coefficients
involved in (\ref{ex}) are rapidly decaying in $k$
and therefore higher order terms can be neglected. In \cite{AbKr} it has been
suggested that $N$ should be of order $h\sqrt{\l}$ where $h$ denotes the mean
width of the duct if the curvature of its boundary is not too large.
In the examples of section \ref{ne} the width of the waveguide varies greatly from point to point. The size of $N$
is determined experimentally and found to depend mainly on the width of the narrowest portion of the waveguide.

Equivalently, we reduce (\ref{sl}) to the Hamiltonian
system of $2N$ equations:
\begin{equation}
J H' = K(x,\l) H, \ x\,\in\Real
\label{hs}
\end{equation}
where
$$
H\,=\,\left( \begin{array}{c}
h \\
Ph' + Q^*h
\end{array} \right) \,\in\,{\Lom}^2(\Real), \ J\,=\,
\left( \begin{array}{cc}
0 & -I\\
I & 0
\end{array} \right),
$$
$$
\ K\,=\,
\left( \begin{array}{cc}
-R \ + \ QP^{-1}Q^* & -QP^{-1}\\
-P^{-1}Q^* & P^{-1}
\end{array} \right).
$$
The system (\ref{hs}) is self-adjoint with $P=P^*>0, \
R=R^*$ for $\l\in\R$. We can therefore apply  advanced
numerical methods (see, for example, \cite{Abr1, Abr2})
to find the \EVs of the problem and the relevant solutions.
Before proceeding to this task let us discuss the issue of boundary conditions.

\section{Boundary conditions}

\subsection{Self-adjoint problem}

To make sure the original \BCs are involved in the ODE
problem consider a generic situation when we have a
functional
$$
J(f) = \int_{\O_0} F(x,y,f,f_x,f_y) \rmd x \rmd y.
$$
To derive the corresponding Euler equation we replace $f$
by $f+\e \g$ and compute
$$
\d J  \ = \ \e \int_{\O_0} \g \left( F_f \  - \ \frac{\dpa}{\dpa x} F_{f_x} \ - \ \frac{\dpa}{\dpa y} F_{f_y}
\right) \rmd x \rmd y
$$
$$
+ \ \e \int_{\dpa \O_0} \g \left( F_{f_x} \rmd y
\ - \ F_{f_y} \rmd x \right).
$$
Putting the first integral equal to zero we obtain the
differential equation (\ref{ee}); the second is
responsible for boundary conditions. For our class
of problems $\dpa \O_0 = \{y=0\} \cup \{y=1\} \cup
\{x=a\}\cup \{x=b\}$. The conditions at different parts of
the boundary are defined by
$$
\int_a^b \left. \g F_{f_y} \right|_{y=0;\,1} \rmd x \ = \ 0;
$$
$$
\int_0^1 \left. \g F_{f_x} \right|_{x=a;\,b} \rmd y \ = \ 0.
$$
Taking account of the obtained quadratic form, we get
\begin{eqnarray}
\int_a^b  \left. \g f_y\right|_{y=0} \rmd x \ = \ 0, \label{bc1} \\
\int_a^b  \left. \g \left(\f'(x) f_x - \frac{1+\f^{\prime 2}}{\f} f_y\right)\right|_{y=1} \rmd x \ = \ 0
\label{bc2}
\end{eqnarray}
on the horizontal lines.
As pointed out in subsection~\ref{sb1}, \DBC remain
unchanged in the new variables and are
automatically taken into account by virtue of our choice
of the functions $g_k(y)$ in (\ref{ex}). The above
integrals (\ref{bc1}), (\ref{bc2}) vanish because of the implied condition $\g=0$.
A difficulty would only occur if we had more
complicated conditions at the curvilinear part of the
boundary of $\O$. \DBC are the ones relevant for quantum mechanical problems and they
enable us to separate the variables in the quadratic form explicitly.
We refer to \cite{AbKr} where the authors consider arbitrary \BCs
of the form $\left(\left. a f + b \frac{\dpa f}{\dpa n}\right)\right|_{\dpa\O}
= 0$ by using appropriate orthogonal curvilinear coordinates. The problem of this kind
requires a more complicated expansion to be used instead of (\ref{ex}). In that case Fourier
coefficients are not obtained in closed form but should be calculated numerically.

On the vertical parts of the boundary we have
\begin{equation}
\int_0^1 g \left. (\f f_x - \f' y f_y)\right|_{x=a,\,b} \rmd y \ = \ 0.
\label{bc3}
\end{equation}
The Dirichlet case is as easy to treat as before: the
conditions $h_k(a) \,=\,h_k(b)\,=\,0$, $k=1,\ldots,N$ are imposed on the solutions
of (\ref{hs}). Consider also a domain where
$\f'(a)\,=\,\f'(b)\,=\,0$ --- the situation typical for
compactly perturbed strips and, in particular, for some waveguides.
Here we are able to handle a more general case. For
instance, \NBC at $x=a,\ b$ do not change and become
$h_k'(a)\,=\,h_k'(b)$, $k=1,\ldots,N$ in terms of the system (\ref{hs}).
However, one cannot fully separate variables in generic Robin conditions
of form (\ref{bc3}).

\subsection{Radiation condition}\label{rc}

Bearing in mind the resonance problem that is of main
interest to us, let us consider a domain $\O$ such that, in notation of subsection \ref{sb1}
$x\in(a,\infty)$ and
\begin{equation}
\f(x)\,\sim\,1, \qquad x\,\geq\, X \label{str}
\end{equation}
for some $X>0$. Similar assumptions are often made in papers dealing with scattering problems,
for instance in \cite{AbKr, OtLa,APV}. We take \DBC on $\partial\O$ and require
a different type of condition to be satisfied as $|x|\to\infty$.
Namely, for a given $\l\in\C$ there always exists a unique solution of (\ref{hel})
that has the form
\begin{equation}
f(x,y) = \left(\exp(-t_1 x) + s_1 \exp(t_1 x)\right) g_1(y) +
\sum\limits_{k=2}^{\infty} s_k g_k(y) \exp(t_k x), \ x\geq X.
\label{rad}
\end{equation}
Here we denote
$$
t_k = -\sqrt{(\pi k)^2-\l}, \qquad \re\, t_k<0;
$$
$g_k$ are the same as in subsection \ref{sb2}.
The coefficients $s_k$, $k=1,2,\ldots$ are defined
by the formula (\ref{rad}) uniquely for each value of $\l$.
We put $\o=\sqrt{\l}$ and consider $s_k$ as functions
of $\o$. The function $s_1$ called the scattering
coefficient of the problem is involved in the definition
of resonances. The reader will find their general definition in
\cite{CFKS}. Note that when variables are separated the
following construction proves to be more handy.
If the scattering coefficient $s_1(\o)$ has a pole at
$\o=\o_0$ we say that $\l=\o_0^2$ is a resonance.
This does not include all the resonances but only those
lying on the first non-physical sheet (see \cite{CFKS}
for detailed explanation).
We refer to \cite{SRMPP,APV} for the equivalence of the
two definitions. Apart from its simplicity, the approach
based on (\ref{rad}) has another distinctive advantage.
It is known from scattering theory that
$s_1(\o)\bar{s_1}(\bar{\o})=1$ and that $s_1$  is analytic
in the half-plane $\im\o<0$. Therefore instead of seeking
the poles of $s_1(\o)$ one can look for its zeros located in the
lower $\o$-half-plane. This is the approach we use
here along with the separation of variables in the perturbed cylinder $\O$.

Given the above definition, there is an obvious
difference between the resonance problem and a classical
spectral problem. Indeed, according to (\ref{rad}) here
we are looking for a solution exponentially growing at
infinity. Note, however, that the resonances we are
interested in occur as perturbations of \EVs and are
typically situated near the real axis. This
means that the values of $|\im\, t_k|$, $k=1,2,\ldots$,
are rather small and therefore the corresponding solution
grows slowly.

Combining (\ref{rad}) with (\ref{ex}) for a sufficiently
large $x$ we get
$$
h_1(x)\,=\,\exp(-t x) \,+\,s_1\exp(t x), \qquad
h'_k(x)\,=\,t_k h_k(x), \ k=2,3,\ldots .
$$
There are two ways to handle these conditions.
One can solve the inhomogeneous problem (as has been done in \cite{APV} in two
dimensions), then find the zeros of $s_1$. Alternatively
one can put $s_1=0$ straight away, then solve the resulting
\EVP with $\l$-dependent boundary conditions. The latter
approach leads to the set of \BCs at $X$
\begin{equation}
\psi_X H(X)\,=\,0, \ \psi_X\,=\,(T-P^{-1}Q^*, \ P^{-1}), \
T = {\rm diag} (t_1,-t_2\ldots,-t_N).
\label{BC}
\end{equation}
This formula known as the radiation, or outgoing wave
condition singles out the solution whose first component grows
and the others decay exponentially at infinity. It is
this solution that is sometimes called the resonance
eigenfunction.

It should be mentioned here that our approach
agrees with that using exterior complex
scaling (see, for example, \cite{CFKS}).
In this technique one replaces the operator by a family of
operators on the same domain which depend analytically on a complex parameter.
The operators are independent of the parameter for $x\leq X$ and are
associated with a space scaling for $x\geq X$. One computes the
complex eigenvalues of this family acting in $\Lom^2(\O)$ and
proves that they do not depend on the parameter, subject to certain conditions.
It is known that the complex eigenvalues coincide with resonances
defined via either the scattering coefficients or analytic continuation
of the resolvent kernels. One may verify that exterior complex
scaling yields the same \BC at $x=X$ as (\ref{BC}).

\section{Numerical examples}\label{ne}

\subsection{Transfer method}\label{tr}
Summarising the results of the first three
sections let us formulate the problems to be solved
numerically. We are looking for such values of $\l$ that the
system (\ref{hs}) has a non-trivial solution satisfying
$$
\psi_a H(a)\,=\,0, \qquad \psi_b H(b)\,=\,0
$$
where
\begin{enumerate}
\item[{\bf 1}] $\psi_a\,=\,\psi_b\,=\,(I, \,0)$;
\item[{\bf 2}] $\psi_a \,=\,(I,\, 0)$, $\psi_b\,=\,\psi_X$ as defined
by (\ref{BC}).
\end{enumerate}
Problem {\bf 1} provides approximations to the Dirichlet \EVs of a compact
domain of kind (\ref{dom}); problem {\bf 2} enables us to calculate complex
resonances that may occur in an unbounded domain of the same
type satisfying (\ref{str}).

The method we apply to both problems is based on the
orthogonal transfer of \cite{Abr1} that we shall
briefly outline below. The manifold of the solutions of the system (\ref{hs})
satisfying the left \BC is determined by
$$
\psi(x) H(x) = 0, \qquad a\leq x \leq b,
$$
where $\psi\in\C^{N \times 2N}$ solves the Cauchy problem
$$
\psi' = \psi J K, \qquad \psi(a) = \psi_a.
$$
Theoretically one can integrate the above equation for a fixed $\l$, define
$$
f(\l) = \det\left(
\begin{array}{c}
\psi(b; \l)\\
\psi_b(\l)
\end{array}
\right)
$$
and solve $f(\l)\,=\,0$ to find the \EVs of the problem.
This method is known to be hopelessly inefficient because
$\psi(x)$, although formally of rank $N$, can have
almost linearly dependent rows. Abramov \cite{Abr1} proposed
replacing $\psi$ by $\tilde{\psi}(x)=\n(x)\psi(x)$, where
$\n\in\C^{N \times N},\,{\det \n}\not=0$. The
function $\n$ is chosen so that to ensure
$\tilde{\psi}(x)\tilde{\psi}^*(x)={\rm const}$.
The transfer equation now takes the form
\begin{equation}
\tilde{\psi}' + \tilde{\psi} JA \left(I -
\tilde{\psi}^*(\tilde{\psi}\tilde{\psi}^*)^{-1}\tilde{\psi}\right)
= 0,  \qquad \tilde{\psi}(a) = \psi_a. \label{ortr}
\end{equation}
The RHS of (\ref{ortr}) is bounded
and the solution $\tilde{\psi}$ exists on the whole of
$[a,b]$. By comparison with  $\psi$, the matrix
$\tilde{\psi}$ has the key advantage of being easily computed
without loss of rank. The use of this idea proved
essential to obtain stable results for this problem.

Having calculated the smooth function $\tilde{\psi}$ we
proceed to find the eigenvalues. For the resonance problem this is done along
the lines of \cite{AYu} where the idea of \cite{Abr1} has
been applied to \nsa \EV problems. As observed there,
$$
\tilde{f}(\l) = \det\left(
\begin{array}{c}
\tilde{\psi}(b)\\
\psi_b
\end{array}
\right) = f(\l) \det\n,
$$
so that the zeros of $f$ and $\tilde{f}$ coincide.
Moreover, the zeros of $\tilde{f}$ can be found by the
method based on the argument principle although $\tilde{f}$
does not have to be analytic in $\l$ as opposed to $f$.
Still the number of zeros of $\tilde{f}$ inside a
contour $\G$ is
$$
N = \frac{1}{2\pi} \oint_{\G} d \Arg \tilde{f}(\l)
$$
as shown in \cite{AYu}. It is this  computational
formula that we use to locate the complex \EVs of problem
{\bf 2}. Taking a shrinking sequence of contours $\G$ we find the zeros up to a chosen accuracy.
We computed the contour integrals reliably for circles
of radii down to $10^{-4}$ and made sure that if the
centre was shifted by a similar order of magnitude the
integrals vanished. In the better conditioned case {\bf 1} we
applied Newton's method allowing us to calculate the \EVs
of the \sa problem.

Note that typically problem {\bf 2} is much harder to solve
than {\bf 1}, and our examples are no exception. When resonances
are situated near the real axis the scattering coefficient
$s_1$ has a pole and a zero close to one another. Naturally, the closer
they are the less stable the problem is.

\subsection{Results of computations}

As an example we consider a quantum waveguide with indentations
defined in the Cartesian coordinates $(\x,\y)$ as
$$
W = \{ -\infty < \x < \infty, \
0 < \y < \f(\x)= 1 - \a \left(\rme^{-(\x-\g)^2} +
\rme^{-(\x+\g)^2}\right)\}
$$
where $\a$ and $\g$ are real positive constants (see fig.~1.1).

We shall be working in two different domains:
$$
\O_{\bf 1}\,=\,W\cap \{0< \x < \g\}, \qquad
\O_{\bf 2}\,=\,W\cap \{0< \x\}.
$$
The domains $\O_{\bf 1}$ and $\O_{\bf 2}$ relate to problems {\bf 1} and
{\bf 2} of the previous subsection, respectively. In the latter the resonance \BC is
imposed at a sufficiently far point $X$ as suggested in
subsection \ref{rc}. In our experiments we put $\g=2$, so that it suffices to take $X \geq 5.4$
to ensure $\f'(X)<10^{-4}$.

In this example there exists $\a=\a_*\approx1$ such that
the two parts of $\dpa\O$ touch one another near $\x=\g$. For this value of
$\a$ the domain $\O$ consists of three disjoint parts
as shown in  fig.~1.2, so that the \EVP is decomposed into three
separate problems. The Laplacian considered
in the compact domain $\O_*$ has
infinitely many real \EVs accumulating at infinity. As
we decrease $\a$ joining the three subdomains, one expects
the \EVs to disappear generating resonances in their neighbourhood.
A similar phenomenon where resonances originate from
\EVs as the domain is perturbed has been observed
in \cite{APV} although the mechanism by which they
emerge is different here.

Clearly, both domains $\O_{\bf 1}$ and $\O_{\bf 2}$ satisfy the
conditions of subsection \ref{sb1} for $\a<\a_*$. We
compute resonances, \ie \EVs of problem {\bf 2} for a range
of $\a$ using the perturbed cylinder approach.
In the \sa example the \EVs of problem {\bf 1} are found for
the same values of $\a$. The ODE problem is solved
by the transfer method described in the previous subsection.
An important question is how to choose the number
of terms to be retained in (\ref{ex}) or, in other words,
the dimension of the system (\ref{hs}). This number
should depend on $\a$: indeed, for $\a$ close to $\a_*$
the width of the domain is small near $\x=\g$
so that one needs to keep a larger number of terms $g_k(y)$ in
(\ref{ex}). There are two possibilities here: first, to
increase repeatedly the number of terms by one and
solve the system with the corresponding constant number of
unknowns until the answers converge. Secondly, instead of keeping
a large number of terms throughout the interval one can
start off with a smaller $N$ in (\ref{hs}).  Moving along the $x$-interval
one changes $N$ gradually adding
or removing variables $h_k(x)$ depending on the size of $\f(x)$.
We have been using both techniques in different
situations ensuring that the results coincide within
the chosen accuracy for two subsequent values of $N$. For reasonably small
values of $\a$ it suffices to take smaller $N$: for
instance, when $\a=0.8$ the results obtained for $N=4$ and
greater coincide up to the tolerance of $10^{-4}$.
The maximal number of terms taken in our computations is
$N=30$ for $\a=0.97$ (the largest value considered).

There is an important connection between the two problems {\bf 1} and {\bf 2} which makes us
study them within the same framework. Namely, as $\a\to\a_*$
both the \EVs of $\O_{\bf 1}$ and the resonances of $\O_{\bf 2}$ converge
to the Dirichlet \EVs of the domain $\O_*$. These
cannot be found by the same method because $\O_*$ has a
cusp at $\x=\g$ and is not a perturbed cylinder in our
terminology. In this case one could still separate variables in the Helmholtz
equation in a similar way arriving at a singular ODE
problem. This question requires a separate consideration
which is beyond the subject of our paper.
Recall that the aim of our numerical experiments
is to calculate resonances occurring in this example
and find out how they are related to eigenvalues.
To be able to make proper comparisons we have used the finite
volume method to discretise the operator on $\O_*$ and find its
spectrum. To know the limit \EV is also helpful as
it  serves as an initial guess for the \EVs of problems {\bf 1}
and {\bf 2}.

We applied the finite volume method to a
series of problems of type {\bf 1} to compare the effectiveness
of the two approaches. As expected, the comparisons are in
favour of the discussed method, which appears to be several
times faster than the conventional one. The benefits of
our technique are more spectacular for larger values of
$\a$. For example, when $\a=0.8$ it takes
twice as long to get accurate results by the finite
volume method, while for $\a=0.9$ the new method
is almost four times faster. Everything else being
equal, the closer $\a$ gets to its critical value the
more advantageous the perturbed cylinder approach is.

The numerical results presented below are related to
the lowest Dirichlet \EV of $\O_*$ and quoted
in terms of the wave number $\o=\sqrt{\l}$
(we shall retain the term \EV for the wave numbers). The smallest
\EV corresponding to the domain with the cusp is
$\o_*=4.6252$; it is included in the diagrams to illustrate
the convergence of our results.

In fig.~2 a series of the \EVs of problem {\bf 1} for
$\a\in[0.7,0.97]$ is shown. They
converge to $\o_*$ as $\a\to\a_*$. The linear rate of convergence is
in agreement with standard perturbation theory: the
principal correction term is of order $\e=\a_*-\a$
since the coefficients of (\ref{hs}) depend on $\a$
linearly.

The resonances $\o_{\bf 2}$ emerging from the lowest \EV \ as $\a$
decreases can be found in table~1. They are also shown in fig.~3 where
their real parts are plotted against their imaginary
parts. They definitely converge to $\o_*$. The problem is, of
course, very sensitive to perturbations near the critical value of
$\a=\a_*$, or $\e=0$. However, we believe that for $\a\leq 0.97$ our
computations are relatively stable and provide reliable results.
For larger values of $\a$ we have observed various
instability effects preventing us from computing resonances
accurately. The system (\ref{hs}) is difficult to solve for values of
$\o$ close to $\o_*$ (and, consequently, the
corresponding resonance) because its coefficients change
very rapidly for $\a\approx\a_*$. Even methods suitable
for stiff systems fail to produce satisfactory results
when $\e$ is too small. As discussed in the previous
subsection, another reason why the problem is likely to
be unstable is the closeness of resonances to the
real axis. These two factors make calculations
slow and inefficient for $\a$ close to $\a_*$.

{\small
\begin{table}
\caption{Resonances $\o_{\bf 2}$ }
\label{t1}
\begin{tabular}{|c|c|c|c|}
\hline
$\a$ & $\o_{\bf 2}$ &$\a$ & $\o_{\bf 2}$ \\\hline
0.7 & $4.2988+0.0545 i$ & 0.9 & $4.5250 +0.0050i$  \\\hline
0.75 & $4.3715 +0.0348i$ & 0.925 & $4.5492+0.0032 i$\\ \hline
0.8 & $4.4223 +0.0212i$ & 0.95  & $4.5755 +0.0017i$ \\ \hline
0.825 & $4.4498 +0.0154i$ & 0.96 & $4.5864+0.0012i$\\ \hline
0.85 & $4.4741+0.0113 i$ & 0.97 & $4.5967 +0.0008 i$\\ \hline
\end{tabular}
\end{table}
}

The $\a$-dependence of the real and
imaginary parts of the resonance originating from $\o_*$
is shown in fig.~4. The real part is found to depend on
the perturbation parameter linearly, whereas the imaginary part seems to be
of order $\e^p$ with $p\approx3/2$ for $\e<0.2$. A different
behaviour has been observed in \cite{APV}, where for a different
geometry the authors derived an asymptotic  formula
for the imaginary part of the resonance. The resonances
in the example of that paper are shown to be analytic
in the  domain perturbation parameter, their imaginary parts depending
quadratically on the parameter.

One generally expects that if one perturbs an \EV which is embedded in the
continuous spectrum, then it is transformed into a resonance near the real
axis. In many cases one can even write out a perturbation expansion, and a
general theory for some such cases was described by Agmon \cite{Ag}. However in
our situation the natural parameter $\e$ can only take positive values
for obvious reasons. Therefore even if the resonance depends analytically on $\e$
for $\e>0$, there is no proof that it has an expansion with finite
coefficients around $\e=0$, nor even that the resonance converges to the \EV as
$\e\downarrow0$. The numerical experiments do, however, suggest that not only does it converge,
but also that its real part and imaginary parts are expanded in powers of $\e^{1/2}$, $\e\downarrow0$.

In the table below the \EVs of problem {\bf 1} are compared with the real parts of the resonances
related to the same values of $\a$. We tabulate the difference $\d=\o_{\bf 1} - \re\,\o_{\bf 2}$ between the
Dirichlet \EV and the real part of the associated resonance and observe that this quantity
decreases as $\e\downarrow0$. One might hope to deduce the rate of convergence of $\d(\e)$
from asymptotic perturbation formulae. The question how to obtain such formulae for the resonance still
remains open.

\begin{table}
\caption{The values of $\d=\o_{\bf 1} - \re\,\o_{\bf 2}$}
\label{t2}
\begin{tabular}{|c|c|c|c|c|c|}
\hline
 $\e$ & 0.3 & 0.25 & 0.2 & 0.175 & 0.15  \\
$\d$ & 0.0436& 0.0256& 0.0182 & 0.0144 & 0.0123\\\hline
 $\e$ & 0.1 & 0.075 & 0.05& 0.04& 0.03 \\
$\d$  &0.0081 & 0.0044 & 0.0015 & 0.0008  & 0.0004\\\hline
\end{tabular}
\end{table}

\section{Conclusions}
The main idea of this paper was  to
reduce an \EVP in a two-dimensional perturbed cylinder to
an ODE problem. Complex resonances occurring in
perturbed strips were also dealt with in the same way.

The method of this paper allowed us to discretise the
problem in one direction taking into account the
geometry of the domain. This was especially important
for the irregularly shaped domains of section \ref{ne}. Standard finite difference
methods would require a significant mesh refinement in
the narrow part of the considered waveguide. For
comparison purposes we also computed the \EVs of
problem {\bf 1} by the finite volume method. As discussed in the
previous section, it proved to be substantially more
time-consuming than the method based on the
separation of variables for the class of problems
studied here. This is in agreement with the already mentioned
results of \cite{AbKr,OtLa} where other advanced
methods were designed to suit similar problems.
We believe that our approach is
competitive and recommend it for ill-conditioned \EV and
resonance problems. Our confidence is supported by the
strong agreement between numerical results and analytic
expectations. Indeed, the convergence of both \EVs and
resonances to their limit value (computed by an independent
method) confirms the reliability of the proposed technique.

\pagebreak

\begin{figure}[h]
\begin{picture}(100,150)(0,0)
\includegraphics{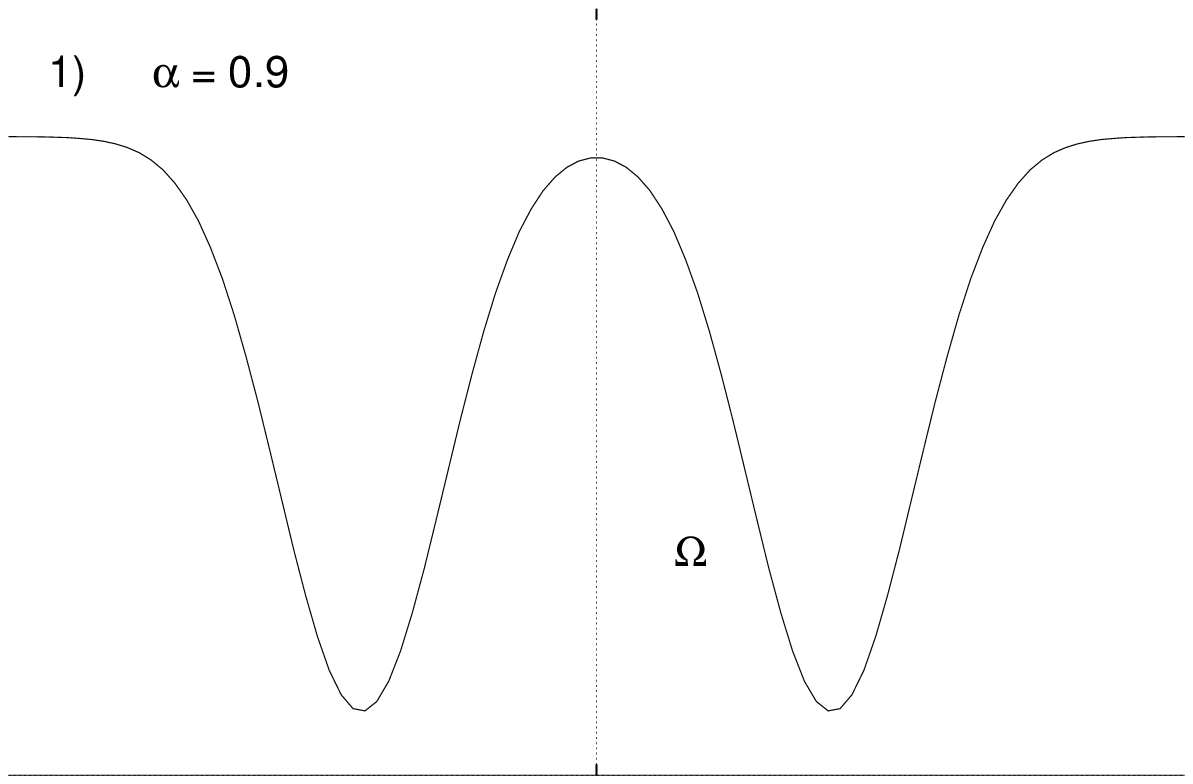}
\includegraphics{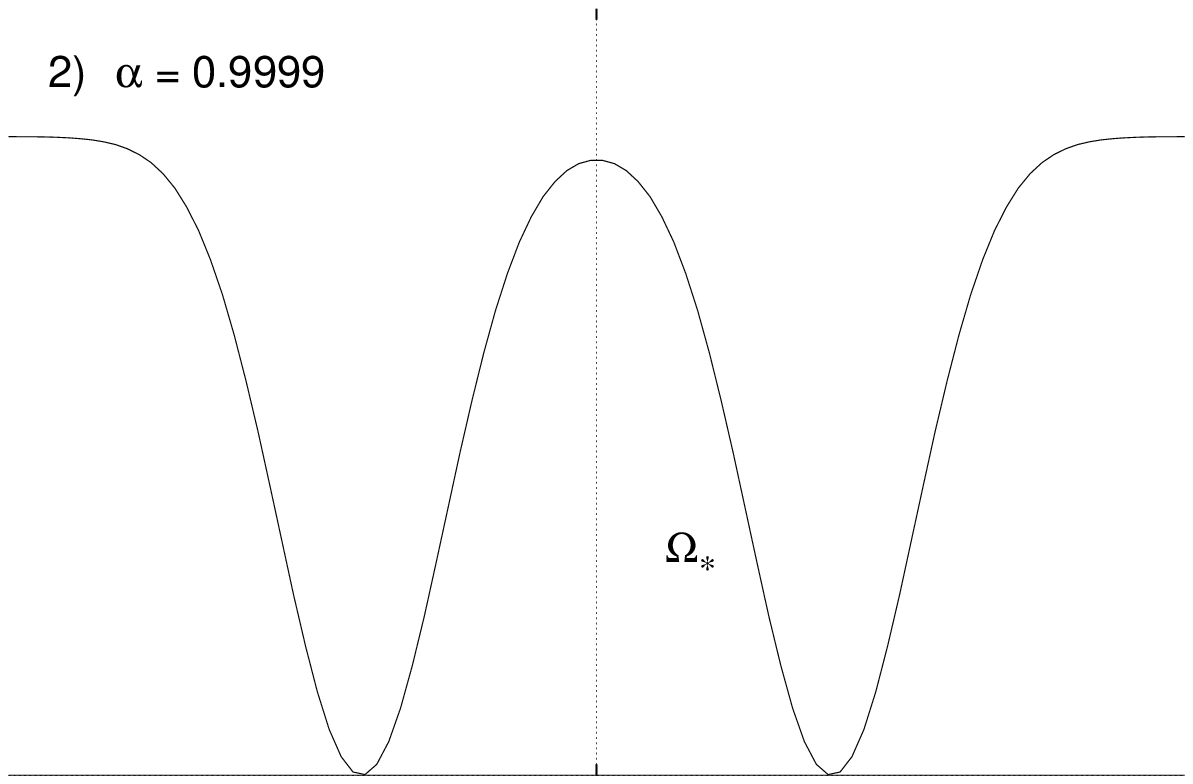}
\end{picture}
\caption{Waveguide with indentations}
\end{figure}


\begin{figure}[h]
\begin{picture}(150,150)(0,0)
\includegraphics{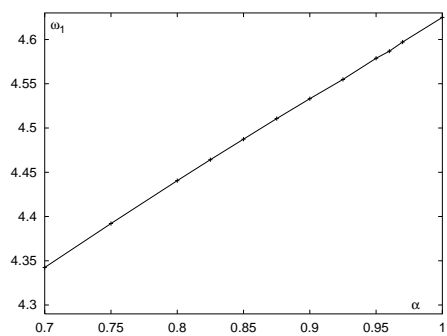}
\end{picture}
\caption{Eigenvalues of problem {\bf 1} for a range of $\a$}
\end{figure}

\begin{figure}[h]
\begin{picture}(150,150)(0,0)
\includegraphics{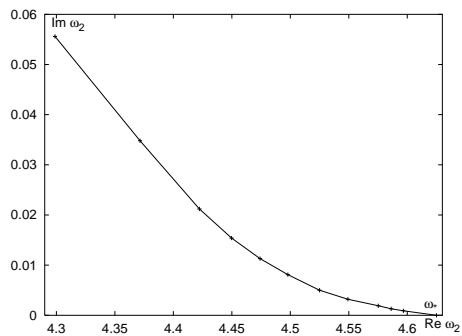}
\end{picture}
\caption{Resonances generated by $\o_*$}
\end{figure}

\pagebreak

\begin{figure}[h]
\begin{picture}(300,150)(0,0)
\includegraphics{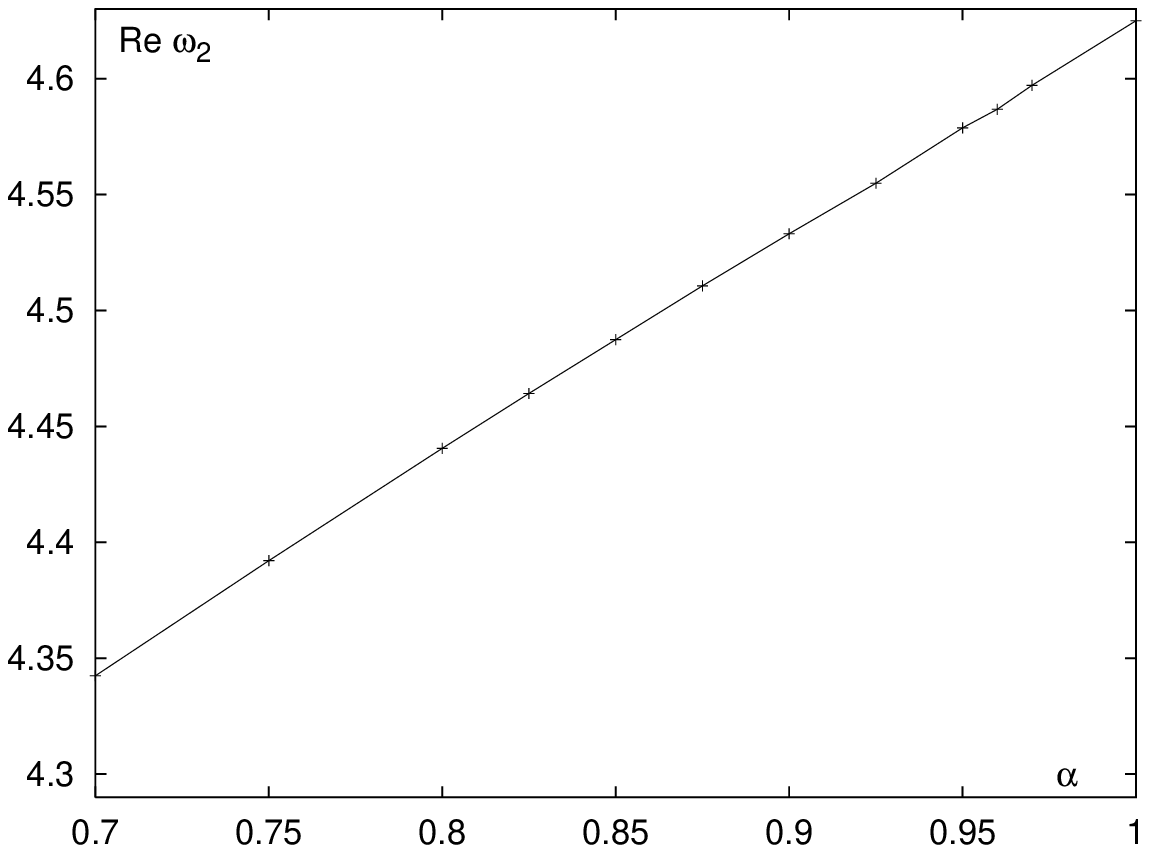}
\includegraphics{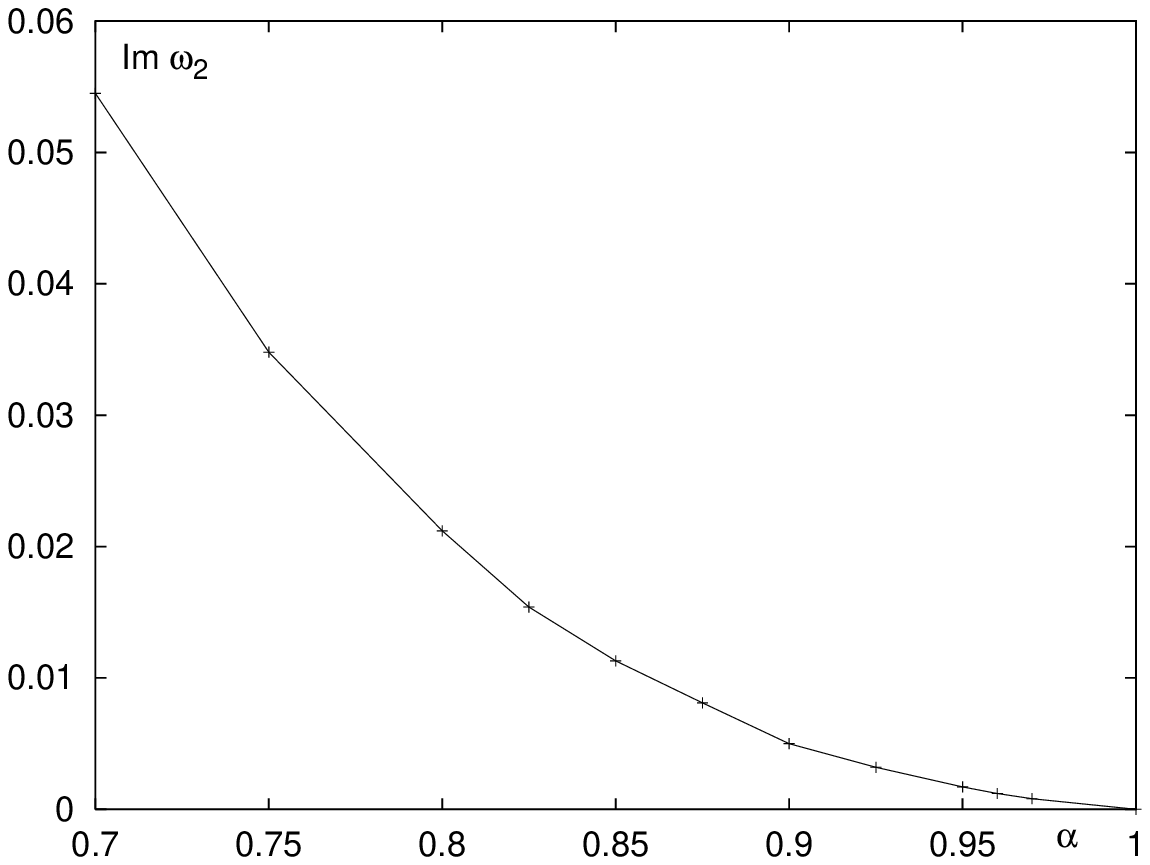}
\end{picture}
\caption{Real and imaginary parts of resonances: $\a$-dependence}
\end{figure}



\end{document}